\documentclass[a4paper,11pt]{article}
\usepackage[T1]{fontenc}
\usepackage[utf8]{inputenc}
\usepackage{authblk}
\title{More than one Author with different Affiliations}
\author[1]{Marco Ferrante\thanks{marco.ferrante@unipd.it}}
\author[2]{Alessia Tagliavini\thanks{alessia.tagliavini@dei.unipd.it}}
\affil[1]{Department of Mathematics, University of Padua, Italy}
\affil[2]{Department of Information Engineering, University of Padua, Italy}
\usepackage[english]{babel}
\usepackage[psamsfonts]{amssymb}
\usepackage{cmmib57}
\usepackage{exscale}
\usepackage{amsmath}
\usepackage{amscd}
\usepackage{latexsym}
\usepackage{graphicx}
\usepackage{amsmath}
\newtheorem{theorem}{Theorem}

\newtheorem{example}[theorem]{Example}

\newcommand{\NN}{\mathbb{N}}

\newcommand{\PP}{\mathbb{P}}
\newcommand{\EE}{\mathbb{E}}

\newcommand{\no}{\noindent}

\def\squarebox#1{\hbox to #1{\hfill\vbox to #1{\vfill}}}

\begin{document}
\title{On the coupon-collector's problem
with several parallel collections}
\maketitle
\begin{abstract}
\no
In this note we evaluate the expectation and variance of the
waiting time to complete $m$ parallel collections of coupons, in the
case of coupons which arrives independently, one by one
and with equal probabilities.
\end{abstract}

\vspace{1truecm}
  {\bf AMS Classification:} 60C05

%\vspace{1truecm}
%  {\bf Short title:} Coupon collector

\section{Introduction}
The coupon-collector's problem is a classical problem in combinatorial probability.
The description of the basic problem is easy:
consider one person that collects coupons and
assume that there is a finite number, say $N$, of different types of coupons.
These items arrive one by one in sequence,
with the type of the successive items being independent random variables
that are each equal to $k$ with probability $p_k$.
%It is immediate to see how this description can be adapted
%to the general problem to draw independent samples
%from a given, finite distribution.

In the coupons-collector's problem,
one is usually interested in answering the following questions:
which is the probability to complete the collection
(or a given subset of the collection)
after the arrival of exactly $n$ coupons ($n \ge N$)?
which is the expectation and the variance of the number of coupons that 
we need to complete the collection (or to complete
a given subset of the collection)? 
which is the expectation and the variance of the number of coupons that 
we need to complete a set of collections?
how these probabilities and expectations change if we
assume that the coupons arrive with unequal probabilities or
in groups of constant size?

The first results, due to De Moivre, Laplace and
Euler (see \cite{MR1082197} for a comprehensive introduction on this topic),
deal with the case of
constant probabilities $p_k\equiv \frac{1}{N}$,
while the first results on the unequal case have to be ascribed to Von Schelling
(see \cite{MR0061772}).
Many other studies have been carry out on this classical problem ever since
(see e.g. \cite{MR0268929}, \cite{MR959649}, \cite{MR1978107} and
\cite{FerSal}).

In this note we will consider the waiting time to complete
a set of $m$ collections, all independent to each other. 
We will assume that any collection is made of a finite number of different coupons,
at any unit of time arrive $m$ new coupons, one for each collection, 
and that the probability to purchase 
any type at any time is uniform. 
We will be able to derive the expectation and variance of the waiting time
to complete all the $m$ collections.

\section{Single collection with equal probabilities}
Let us start by considering a single collection.
Assume that this collection consists of $N$ different coupons, which 
are equally likely, with the probability to purchase 
any type at any time equal to $\frac{1}{N}$. 
In this section we will evaluate the expectation and the variance
of the random number of coupons that 
one needs to purchase in order to complete the
collection. 
We will follow two approaches, the first one present in most of the textbooks in Probability,
while the second one, based on a Markov Chain
approach, will allow us to extend the computation to the case of 
parallel collections. 

\subsection{The Geometric Distribution approach}

Let $X$ denote the (random) number of coupons that we need to purchase in order to complete our collection.
We can write $X=X_1+X_2+\ldots +X_N$, where for any $i=1,2, \ldots, N$, $X_i$ denotes 
the additional number of coupons
that we need to purchase to pass from $i-1$ to $i$ different types of coupons in our collection.
Trivially $X_1=1$ and, since we are considering the case of a uniform distribution, it follows that when $i$ distinct types 
of coupons have been collected, a new coupon purchased will be of a distinct type with probability equal to
$\frac{N-i}{N}$. 
By the independence assumption, we get that the random variable $X_i$, for $i \in \{ 2, \dots, N \}$,
is independent from the other variables and has a geometric law with parameter $\frac{N-i+1}{N}$.
The expected number of coupons that we have to buy to complete the collection will be therefore
\begin{align*}
\mathbb{E}[X]	& = \mathbb{E}[X_1] + \dots + \mathbb{E}[X_N] 
%			= 1 + \frac{N}{N-1} + \frac{N}{N-2} + \cdots + \frac{N}{2} + N \\
%			& 
			= N \sum_{i=1}^N \frac{1}{i}
\end{align*}
while its variance
\begin{align}
\label{varExp}
Var [X]	
%& = Var [X_1] + \dots + Var [X_N] 
		 = \sum_{i=0}^{N-1} \frac{N i}{(N-i)^2}
%		 \nonumber \\
%			& 
			= N^2 \sum_{i=1}^N \frac{1}{i^2} - N \sum_{i=1}^N \frac{1}{i}
\end{align}

\subsection{The Markov Chain approach}
Even if the previous result is very simple and the formula completely clear, 
it will not help us in order to deal with the problem of parallel collections.
We will introduce then a different
approach, that will provide an alternative to the previous computation and
that we will be able to extend in the following section, where the situation becomes more
complicate.

Assuming as before that one coupon arrives at any unit of time,
it is possible to solve the previous problem by using a Markov Chains approach.
Indeed, let $Y_n$ be the number of different types of coupons collected after $n$ units of time and 
assume again that the probability of finding a coupon of any type at any time is $p = \frac{1}{N}$. 
$\{Y_n,n\in \NN\}$ is then a Markov Chain on the state space $S = \{ 0, 1, \dots, N \}$ 
%with $|S|=N+1$
and it is immediate to see that its
transition matrix is given by
$$
P = 
\begin{pmatrix}
0 & 1 & 0 & \dots & \dots & \dots & 0 \\
0 & \frac{1}{N} & \frac{N-1}{N} & 0 & \dots & \dots & 0 \\
0 & 0 & \frac{2}{N} & \frac{N-2}{N} & 0 & \dots & 0 \\
\vdots & \vdots & \ddots & \ddots & \ddots & \ddots & \vdots \\ 
0 & \dots & \dots & \dots & 0 & \frac{N-1}{N} & \frac{1}{N} \\
0 & \dots & \dots & \dots & 0 & 0 & 1
\end{pmatrix}
$$
Note that $\{N\}$ is the unique closed class and the state $N$ is therefore absorbing,
while the other states are transient. If we denote by $Q$ the sub matrix of $P$
relative to the transient states (in this case the first $N$ lines and rows) 
we can define the fundamental matrix of $\{Y_n,n\in \NN\}$ as
\[
F=(Id-Q)^{-1}
\]
Note that $QF=FQ=F-Id$ (see \cite{MR0410929} for more details).

If we define the random variable
\[
D^N=\inf\{n\ge0:Y_n=N\}
\]
its conditional expectation $\EE_0[D^N]=\EE[D^N|Y_0=0]$ will be equal to the expected number of 
coupons needed to complete the collection.
If we define the vector $k^N=(k^N_0, \ldots, k^N_N)=(\EE_0[D^N],\ldots,\EE_N[D^N])$,
a classical result (see \cite{MR1600720})
states that $k^N$ is the minimal non negative solution of the linear system:
$$
\begin{cases}
k^N_N = 0 \\
k_i^N = 1 + \sum_{j \ne N} p_{ij} k_j^N, \quad i \ne N
\end{cases}
$$
Therefore
\[
(k_0^N,\ldots,k_{N-1}^N)=F 
\left(
\begin{array}{c}
1 \\ \vdots \\ 1
\end{array}
\right)
\]
%In this simple case, we can solve explicitly this system 
%$$
%\begin{cases}
%k_0 = k_1 + 1 \\
%k_1 = \frac{1}{N} k_1 + \frac{N-1}{N} k_2 + 1 \\
%k_2 = \frac{2}{N} k_2 + \frac{N-2}{N} k_3 + 1 \\
%\vdots \\
%k_{N-2} = \frac{N-2}{N} k_{N-2} + \frac{2}{N} k_{N-1} + 1 \\
%k_{N-1} = \frac{N-1}{N} k_{N-1} + 1 \\
%k_N = 0 
%\end{cases}
%\quad \iff \quad
%\begin{cases}
%k_N = 0 \\
%\frac{1}{N} k_{N-1} = 1 \\
%\frac{2}{N} k_{N-2} = \frac{2}{N} k_{N-1} + 1 \\
%\vdots \\
%\frac{N-2}{N} k_2 = \frac{N-2}{N} k_3 +1 \\
%\frac{N-1}{N} k_1 = \frac{N-1}{N} k_2 +1 \\
%k_0 = k_1 + 1
%\end{cases}
%$$
%$$
%\iff
%\begin{cases}
%k_N = 0 \\
%k_{N-1} = N \\
%k_{N-2} = k_{N-1} + \frac{N}{2} = N + \frac{N}{2} \\
%k_{N-3} = k_{N_2} + \frac{N}{3} = N + \frac{N}{2} + \frac{N}{3} \\
%\vdots \\
%k_2 = k_3 + \frac{N}{N-2} = N + \frac{N}{2} + \dots + \frac{N}{N-3} + \frac{N}{N-2} \\
%k_1 = k_2 + \frac{N}{N-1} = N + \frac{N}{2} + \dots + \frac{N}{N-2} + \frac{N}{N-1} \\
%k_0 = k_1 + 1 = N + \frac{N}{2} + \dots + \frac{N}{N-2} + \frac{N}{N-1} + 1
%\end{cases}
%$$
It is immediate to see that in this case
\[
F = 
\begin{pmatrix}
1 & \frac{N}{N-1} & \frac{N}{N-2} &  \frac{N}{N-3} & \dots & \dots & N \\
0 & \frac{N}{N-1} & \frac{N}{N-2} &  \frac{N}{N-3} & \dots & \dots & N \\
0 & 0 & \frac{N}{N-2} & \frac{N}{N-3} & \ldots & \dots & N \\
0 & 0 &0 & \frac{N}{N-3} & \dots & \dots & N \\
\vdots & \vdots & \ddots & \ddots & \ddots & \ddots & \vdots \\ 
0 & \dots & \dots & \dots & 0 & \frac{N}{N-1} & N \\
0 & \dots & \dots & \dots & 0 & 0 & N
\end{pmatrix}
\]
Then, for $j=0,\ldots,N-1$,
\[
k^N_j=N \sum_{i=1}^{N-j} \frac{1}{i}
\]
and the expected waiting time to collect all the different coupons is given by
$$
k^N_0 = N \sum_{i=1}^N \frac{1}{i}
$$

The computation of the conditional variance is similar (see again \cite{MR0410929}).
Denoting by 
\[
v^N_i=Var_i[D^N]=\EE[(D^N)^2|Y_0=i]-
(\EE[(D^N)^2|Y_0=i])^2
\]
we have
\[
v^N=
(2F-Id)
\left(
\begin{array}{c}
k^N_0 \\ \vdots \\ k^N_{N-1}
\end{array}
\right)
- 
((k^N_0)^2,\ldots,(k^N_{N-1})^2)
\]
A simple computation gives that
the variance of the waiting time to 
collect all the different coupons is given by
\begin{align}
\label{varMarkov}
v^N_0=N \sum_{i=1}^N \frac{1}{i}-\left(N \sum_{i=1}^N \frac{1}{i}\right)^2
+\sum_{k=1}^{N-1}\left[
\frac{2N^2}{N-k} \sum_{i=1}^{N-k} \frac{1}{i}
\right]
\end{align}
It is possible to see that (\ref{varExp}) and (\ref{varMarkov}) 
are indeed the same.

\section{Parallel collections with equal probabilities}
Let us now consider the case that $m$ different collections are available.
Assume that the $i$-th collection is made of $N_i$ different coupons 
and that the probability to purchase 
any type at any time is equal to $\frac{1}{N_i}$. 
Moreover, let us assume that we purchase simultaneously one coupon of every collection.
In this section we will derive the expectation and the
variance of the waiting time to complete 
all the $m$ different collections. 
%Note that this number will be close to $f(\sum_{j=1}^m N_j)$. 

With the notation of the previous section, let $X^i$ denote be the random number of
coupons needed to complete the $i$-th collection.
The random number of coupons needed to complete all the $m$ collections
will be therefore $\max(X^1,\ldots,X^m)$.
It looks not simple to determine the law of this random variable, or even just to provide a direct
computation of its expectation.
On the converse, adapt to this situation the Markov Chain's approach is quite
simple.

Let $Y_n=(Y^1_n,\ldots,Y^m_n)$ be the number of different types of coupons 
of any collection, collected after $n$ units of time. 
$\{Y_n,n\in \NN\}$ is again a Markov Chain on the state space 
\[
S = \{(0,\ldots,0)\}\cup\{(i_1,\ldots,i_m): i_j\in\{1,\ldots,N_j\}, j\in\{1,\ldots,m\} \}
\]
We have $|S|=1+\prod_{j=1}^m N_j$.
Defined $(\alpha_1,\ldots,\alpha_m)\in \{0,1\}^m$
and assuming that $(i_1,\ldots,i_m)\in S$ and
$(i_1+\alpha_1,\ldots,i_m+\alpha_m)\in S$, we have
that the transition matrix is defined as
\[
\PP[Y_{t+1}=(i_1+\alpha_1,\ldots,i_m+\alpha_m)|
Y_{t+1}=(i_1,\ldots,i_m)]
=
\]
\[
=\prod_{j=1}^m
\left(1-\frac{i_j}{N_j}\right)^{\alpha_j}
\left(\frac{i_j}{N_j}\right)^{1-\alpha_j}
\]
(Note that $\PP[Y_{t+1}=(1,\ldots,1)|
Y_{t}=(0,\ldots,0)]
= 1$).
The transition matrix has again a unique absorbing state
$\{(N_1,\ldots,N_m)\}$. If we denote by $Q$ the sub matrix of $P$
relative to the transient states
we can define the fundamental matrix $F$ of $\{Y_n,n\in \NN\}$ as before
%\[
%F=(Id-Q)^{-1}
%\]
and compute the conditional expectations 
\[
k^{\{(N_1,\ldots,N_m)\}}=\Big(\EE_{(0,\ldots,0)}[D^{\{(N_1,\ldots,N_m)\}}],\ldots,\EE_{(N_1,\ldots,N_m-1)}
[D^{\{(N_1,\ldots,N_m)\}}]\Big)
\]
as
\[
k^{\{(N_1,\ldots,N_m)\}} = 
F 
\left(
\begin{array}{c}
1 \\ \vdots \\ 1
\end{array}
\right)
\]
\begin{example}
Consider three collections, each with six different coupon to be collected.
With the previous notation we get $m=3$, $N_1=N_2=N_3=6$ and
the matrix $F\in[0,1]^{215\times215}$, even if at most eight entries of
every row are not zeros. 
Using Matlab\textsuperscript{\textregistered}  
we are able to compute $F$ in this case and 
the expectation of the waiting time to complete all the three collections is
equal to $20.01$.
%The expected total number of coupons purcheased will be $60.03$.
%Note that 
%the expected number of coupons needed to complete
%a single collection with 6 different coupons is equal to $14.7$ and one
%with 18 different coupons to $62.91$.
\end{example}
In order to compute the variance of the waiting time
to complete all the parallel collections,
we shall operate as in the previous section.
Denoting, for $i\in S\setminus \{(N_1,\ldots,N_m)\}$,
\begin{align*}
v^N_i & = Var_i[D^{\{(N_1,\ldots,N_m)\}}]
\\ 
& = \EE[(D^{\{(N_1,\ldots,N_m)\}})^2|Y_0=i]-
(\EE[(D^{\{(N_1,\ldots,N_m)\}})^2|Y_0=i])^2
\end{align*}
the conditional variances, 
we have
\begin{align*}
v^N = 
(2F-Id)
\left(
\begin{array}{c}
k^N_{(0,\ldots,0)} \\ \vdots \\ k^N_{(N_1,\ldots,N_m-1)}
\end{array}
\right)
- 
% \Big(\EE_{(0,\ldots,0)}[D^{\{(N_1,\ldots,N_m)\}}],\ldots,\EE_{(N_1,\ldots,N_m-1)}
% [D^{\{(N_1,\ldots,N_m)\}}]\Big)
((k^N_{(0,\ldots,0)})^2,\ldots,(k^N_{(N_1,\ldots,N_m-1)})^2)
\end{align*}

\begin{example}(Continue)
Considering again three collections, each with six different coupon to be collected,
the variance of the waiting time to complete all the three collections is
equal to $44.8975$.
\end{example}

%\lstinputlisting{code.m}

%
%\bibliography{coupon} % Here the name of the file .bib
%                       % (bibliography database)
%                        % that must be used.
%\bibliographystyle{plain}

\end{document}